\documentclass[twoside]{article}
\usepackage{amssymb}

\begin{document}

\baselineskip=16pt

\title{\bf The Circular Chromatic Number of the Mycielski's graph $M^t(K_n)$
\thanks{School of Mathematical Sciences
Laboratory of Mathematics and Complex Systems Beijing Normal
University, Beijing 100875, China ;Project 10271017 supported by
NNSFC}}

\date{}

\author{\rm
Zuqiang Ma \ \ \ Junliang Cai\\
{\it School of Mathematical Sciences, Beijing Normal University}
\\
{\it Beijing\ \ 100875,\ \ P.R.China}
\\
({ mazuqiang@yahoo.com.cn ;\ \ \ caijunliang@bnu.edu.cn})}

\maketitle

\centerline{\bf Abstract}

As a natural generalization of chromatic number of a graph, the
circular chromatic number of graphs (or the star chromatic number)
was introduced by A.Vince in 1988. Let $M^t(G)$ denote the $t$th
iterated Mycielski graph of $G$. It was conjectured by Chang, Huang
and Zhu(Discrete mathematics,205(1999), 23-37) that for all $n \ge
t+2,\ \chi_c(M^t(K_n))=\chi(M^t(K_n))=n+t.$ In 2004, D.D.F. Liu
proved the conjecture when $t\ge 2$, $n\ge 2^{t-1}+2t-2$. In this
paper,we show that the result can be strengthened to the following:
if $t\ge 4$, $n\ge \frac{11}{12}2^{t-1}+2t+\frac{1}{3}$, then
$\chi_c(M^t(K_n))=\chi(M^t(K_n))$.

{\bf Keyword:}   circular chromatic number, complete graph,
Mycielski graph

 {\bf MSC(2000):}   05A30

\section{{Introduction}}
The circular chromatic number of a graph is a natural generalization
of chromatic number $\chi(G)$ of a graph,\ introduced and studied by
A.Vince in 1988, as the `star chromatic number'.

{\bf Definition 1$^{[13]}$} Let $k$ and $d$ be positive integers
such that $k \ge 2d$,\ a $(k,d)$-coloring of a graph $G$ is a
mapping
$$
f:\ \ V(G)\longrightarrow \mathbb{Z}_k=\{0,1,\cdots,k-1\}
$$
such that for any $uv \in E(G),d \le |f(u)-f(v)| \le k-d$. {\it The
circular chromatic number} $\chi_c(G)$ of $G$ is defined as
$$
\chi_c(G)=\inf\{\frac{k}{d}\ :\ \mbox{there exists a}\ (k,d)\mbox
{-}\mbox{coloring of}\ G\}.
$$

Obviously, \ a $(k,1)$-coloring of a graph $G$ is just an ordinary
$k$-coloring of $G$. The following properties of circular chromatic
number can be found in [2,13,14]:

{\bf Property 1}$^{[13]}$\ \ $\chi(G)-1 < \chi_c(G) \le \chi(G).$

{\bf Property 2}$^{[13]}$\ \ If $H$ is a subgraph of $G$,\ then
$\chi_c(H) \le \chi_c(G)$.

{\bf Property 3}$^{[2]}$\ \ \ If there is a homomorphism from $G$ to
$H$,then $\chi_c(G) \le \chi_c(H)$.

{\bf Property 4}$^{[14]}$\ \ For any graph $G$,\ $\chi_c(G)$ is a
rational number.\ In fact,\ we have
$$\chi_c(G)=\min\{\frac{k}{d}\ :\ \mbox{there exists a}\
(k,d)\mbox{-}\mbox{coloring of}\ G, k \le |V(G)|, d \le
\alpha(G)\}.$$ where $\alpha(G)$ is the independence number of $G$.

The question of determining a graph $G$ satisfying
$\chi_c(G)=\chi(G)$ was asked by Vince in [13].\ Unfortunately it is
hard to determine whether a given graph $G$ has $\chi_c(G)=\chi(G)$
or not [5]. In spite of this difficulty,\ many classes of graphs
satisfying $\chi_c(G)=\chi(G)$ have been found,\ and the Mycielski
graph is one of the important classes[3,4,10,11,12].

{\bf Definition 2}\quad Let $G$ be a simple graph with vertex set
$V(G)=\{x_1,x_2,\cdots,x_n\}$ and edge set $E(G)$,\ the {\it
Mycielski graph} of $G$,\ denoted by $M(G)$,\ is the graph with
vertex set
$$
V(M(G))=\{x_1,x_2,\cdots,x_n\} \cup
\{x^{'}_1,x^{'}_2,\cdots,x^{'}_n\} \cup \{u\},
$$
and edge set
$$
E(M(G))=E(G) \cup \{x^{'}_ix_j:x_ix_j \in E(G),1 \le i,j \le n\}
\cup \{ux^{'}_i:1 \le i \le n\}.
$$

For each $x_i \in V(G),x^{'}_i$ is called the $twin$ of $x_i$$(x_i$
is also the $twin$ of $x^{'}_i)$,\ and the new vertex $u$ is called
the $root$ of $M(G)$. Let $M^0(G)=G$,\ and define the $t$th
Mycielski graph of $G$ is $M^t(G)=M(M^{t-1}(G))$.

Here, we consider $M^i(G)$ as a subgraph of $M^j(G)=M^{j-i}(M^i(G))$
if $i \le j$.

We know that for any nonempty graph $G,\ \chi(M(G))=\chi(G)+1$.\
Unfortunately,\ no simple characterization of graphs with
$\chi_c(M(G))=\chi(M(G))$ have been found,\ though some graphs that
satisfy this condition were studied in [4,8].\ An important class of
graphs is the Mycielski graph of complete graph $K_n$.\ The
following conjecture was introduced by Chang [3] and confirmed for
the cases of $t=1,2$.

{\bf Conjecture$^{[3]}$} If $t \ge 1, n \ge t+2,$ then
$\chi_c(M^t(K_n))=\chi(M^t(K_n))=n+t.$

For $t \ge 2$ and $n \ge 2^{t-1}+2t-2$,\ the conjecture is verified
[10].\ Recently,\ Simonyi and Tardos proved the conjecture holds
when $t+n$ is even.

In this paper,\ we show that the result in [10] can be strengthened
to the following:

{\bf Theorem 1}\quad If $t\ge 4, n\ge
\frac{11}{12}2^{t-1}+2t+\frac{1}{3}$, then
$\chi_c(M^t(K_n))=\chi(M^t(K_n))=n+t$.

\section{Some definitions and lemmas}
Let $k \ge 2d$,\ Fan introduced $(k,d)$-partition of $G$ [4]:\ a
$(k,d)$-partition of $G$ is a partition $(X_0,X_1,\cdots,X_{k-1})$
of $V(G)$ such that $V(G)=\bigcup\limits_ {i=0}^{k-1}X_i,$\ and for
each $j \in \mathbb{Z}_{k-1},$
$$ X_j \cup X_{j+1} \cup \cdots \cup X_{j+d-1}$$
is an independent set in $G$ (Here it is allowed that
$X_i=\emptyset$). It is easy to see that a $(k,d)$-partition of $G$
is equivalent to a $(k,d)$-coloring of $G$ (Here we call $X_i$
$color\ class$ of $i, 0 \le i \le n-1$). Thus, we
have$$\chi_c(G)=\min\{\frac{k}{d}\ :\  \mbox{there exists a}\ (k,d)\
\mbox{-partition of}\ G\}.$$

Combining with the conclusions in [4] and [6], the following result
holds:

{\bf Lemma 1$^{[4,6]}$}\ Let\
$\chi_c(G)=\frac{k}{d},(k,d)=1$,$(X_0,X_1,\cdots,X_{k-1})$ is a
$(k,d)$-partition of $G$, then $\forall \ i \in \mathbb{Z}_k$, $X_i
\neq \phi $, and $N(X_i) \cap X_{i+d} \neq \phi$.

To work with such complicated graphs as $M^t(G)$, we need to take a
system to name the vertices of $M^t(G)$ .

{\bf Definition 3$^{[11]}$} \ Suppose $x \in V(G)$ , and $t$ is a
positive integer,

1.\ As a vertex of $M^{t-1}(G)$,\ the twin of $x$ in $M^t(G)$ is
called the $t$th twin of $x$,\ denoted by $x^t$,\ $x$ is called {\it
initial vertex} of $M^t(G)$;

2.\ If $i$ and $j$ are positive integers,\ the $i$th twin of $x$ in
$M^i(G)$ is $x^i$,\ the $j$th of $x^i$ in $M^j(M^i(G))=M^{i+j}(G)$
is denoted by $(x^i)^j$,\ which can be simplified as $x^{ij}$ if
there is no ambiguity.

If $i_1,i_2,\cdots,i_n \in \mathbb{N}$,\ $x^{i_1i_2\cdots i_n}$ is
defined recursively by $(x^{i_1i_2\cdots i_{n-1}})^{i_n}$.\ For any
positive integers $i_1,i_2,\cdots,i_n $ such that
$i_1+i_2+\cdots+i_n \le t$,\ $x^{i_1i_2\cdots i_n}$ is called a
$derived\ vertex$ of $x$ in $M^t(G)$(Here, the $i$th twin of $x$ is
the derived\ vertex of $x$ too). The set of derived vertices of $x$
is denoted by $T(x)$:
$$T(x)=\{x^{i_1i_2\cdots i_n}\ :\ i_1>0, \cdots, i_n>0,
i_1+i_2+\cdots+i_n \le t\}.
$$

3.\ For $0<i \le t,$ the new root in $M^i(G)$ formed from
$M^{i-1}(G)$ to $M^i(G)$ is called the $i$th root of $G$,\ denoted
by $u_i$,\ The set of roots and their derived vertices in $M^t(G)$
is denoted by $R(M^t(G))$:$$ R(M^t(G))=\{u_i,u^{i_1i_2\cdots
i_n}_i:0< i \le t,i_1>0,\cdots,i_n>0,i_1+i_2+\cdots+i_n \le t-i\}.
$$

Fig.1 shows $M^2(G).$

For convenience,\ any vertex in $R(M^t(G))$ is called the $root$ of
$M^t(G)$.

{\bf Definition 4}\quad We call $u^{i_1i_2\cdots i_n}_i$ the $s$th
$root$ of $M^t(G)$ if $i+i_1+i_2+\cdots+i_n=s$, and denote
$R_s(M^t(G))$ as the set of all $s$th $roots$ of $M^t(G)$.

By the definition of Mycielski graph, there is a bijection $h$ from
$R(M^{t-1}(G))$ to $R_t(M^t(G))-\{u_t\}$ as follows
$$
h: R(M^{t-1}(G)) \longrightarrow R_t(M^t(G))-\{u_t\}.
$$
For any $u \in R(M^{t-1}(G)), h(u)$ is the $twin$ of $u$ in
$M^t(G)=M(M^{t-1}(G))$. In other words, $h(u^{i_1i_2\cdots
i_{n-1}}_i)= u^{i_1i_2\cdots i_n}_i$, $i_n \neq 0$,
$i_1+i_2+\cdots+i_n=t$.

 \setlength{\unitlength}{8mm}
\begin{picture}(19,6.5)
\put(3.5,3.15){\oval(4,1)} \put(2,3){$\{x_i^2:x_i \in V(G)\}$}
\put(3.5,5.15){\oval(4,1)} \put(2,5){$\{x_i:x_i \in V(G)\}$}

\put(9.5,3.15){\oval(4,1)} \put(8,3){$\{x_i^{11}:x_i \in V(G)\}$}
\put(9.5,5.15){\oval(4,1)} \put(8,5){$\{x_i^1:x_i \in V(G)\}$}

\put(14,3.15){\circle{1}} \put(14,5.15){\circle{1}}
\put(8,1.15){\circle{1}}

\put(5.5,5.15){\line(1,0){2}} \put(11.5,5.15){\line(1,0){2}}
\put(3.5,3.65){\line(0,1){1}}

\put(5.4,4.8){\line(2,-1){2.4}} \put(5.4,3.5){\line(2,1){2.4}}
\put(11.4,4.8){\line(2,-1){2.4}} \put(11.4,3.5){\line(2,1){2.4}}

\put(5.4,2.8){\line(2,-1){2.38}} \put(9.5,2.65){\line(-1,-1){1.13}}
\put(13.9,2.63){\line(-4,-1){5.4}}

\put(13.8,5){$u_1$} \put(13.8,3){$u_1^1$} \put(7.8,1.1){$u_2$}

\put(7.5,0){Fig.1}
\end{picture}

We say $x \in V(G)$ is a universal vertex if for each $y \in
V(G)-\{x\},xy \in E(G)$.\ Concerning the circular chromatic number
of $M^t(K_n)$, the following results holds:

{\bf Lemma 2$^{[6]}$}\ Let $G$ be a graph with $n$ universal
vertices, $n \ge 2$. If $\chi_c(M^t(G))=\frac{k}{d}$, $(k,d)=1$,
then $ (n-3)(d-1) \le 2^t-2$.

{\bf Lemma 3$^{[11]}$}\ Let $G$ be a simple graph,
$V(G)=\{x_1,x_2,\cdots,x_n\}$. The vertex set of $M(G)$ is
$V(M(G))=\{x_1,x_2,\cdots,x_n\}\cup
\{x^{'}_1,x^{'}_2,\cdots,x^{'}_n\}\cup \{u\}$.If $V(G)$ has
$(k,d)$-partitions, then there is a $(k,d)$-partition of $V(G)$ such
that

1.\ $u \in X_0;$ \ \ \

2.\ For some $i$ ,\ $d \le i \le k-d$,$\forall \ x \in V(G)$,\ we
have $\{x,x^{'}\}\subseteq X_i$ if $x \in X_i.$


{\bf Definition 5$^{[11]}$} Let $(X_0,X_1,\cdots,X_{k-1})$ be a
$(k,d)$-partition of $G$. For any $x \in V(G),$ there exist $j\ ,0
\le j \le k-1$ such that $x \in X_j$. The following set
$$
\delta(x)=X_{j-d+1}\cup X_{j-d+2}\cup \cdots \cup X_j \cup X_{j+1}\cup \cdots \cup X_{j+d-1}
$$
is called the $d$-$field$ of $x$ in partition
$(X_0,X_1,\cdots,X_{k-1})$.

It is easy to see that for any $x \in V(G)$, $\delta(x) \cap
N(x)=\emptyset$,\ where $N(x)$ is the adjacency set of $x$. In
general, we use $C(x)$ to denote the color class containing $x$.

Suppose $\{x_i,x_j\}\subseteq V(K_n),i \neq j$,\ by the definition
of $M^t(K_n)$, it is obviously that $T(x_i) \subseteq N(x_j)$ .
Hence, $X_j \subseteq R(M^t(G))$ if $X_j \subseteq \delta(x_i) \cap
\delta(x_{i+1}),$\ because $\forall \ v \in T(x_l),1 \le l \le n,v
\notin X_j$.

{\bf Definition 6}\quad Let $t \ge 1$,$F^{\circ}_t$ is a digraph
with vertex set
$$
V(F^{\circ}_t)=R_t(M^t(G))-\{u_t\},
$$
and arc set
$$
A(F^{\circ}_t)=\{(u^{i_1i_2\cdots i_n}_i,u^{i_1i_2\cdots
i_{n-1}j(i_n-j)}_i): 1 \le i \le t-1,i+\sum_{l=1}^n{i_l}=t,\ i_n \ge
2,\ 1 \le j \le i_n-1 \}.
$$

For convenience of state and proof, we define a directed graph $F_t$
which is the isomorphic graph of $F^{\circ}_t$, by replacing
$u^{i_1i_2\cdots i_n}_i$ as $u^{i_1i_2\cdots i_{n-1}}_i$,\
$i+i_1+i_2+\cdots+i_n=t$(see Fig.2).

\setlength{\unitlength}{8mm}
\begin{picture}(11,11)
\put(2,.5){ \put(0,9.5){$u_1$} \put(2,9.5){$u^1_1$}
\put(3.9,9.5){$u^{11}_1$} \put(0.6,9.6){$\vector(1,0){1}$}
\put(2.6,9.6){$\vector(1,0){1}$} \put(4.7,9.6){$\vector(1,0){1}$}
\put(6,9.5){$\cdots$} \put(6.7,9.6){$\vector(1,0){1}$}
\put(8,9.5){$u^{\overbrace{\scriptstyle{11 \cdots 1}}^{t-3}}_1$}
\put(9.1,9.6){$\vector(1,0){1}$} \put(0.6,9.3){$\vector(3,-4){1.2}$}
\put(10.4,9.5){$u^{\overbrace{\scriptstyle{11 \cdots 1}}^{t-2}}_1$}
\put(2,7.5){$u^2_1$} \put(2.5,7.2){$\vector(2,-1){1}$}
\put(3.7,6.6){$\cdots$} \put(3.9,7.5){$u^{21}_1$}
\put(2.6,7.6){$\vector(1,0){1}$} \put(4.7,7.6){$\vector(1,0){1}$}
\put(6,7.5){$\cdots$} \put(6.7,7.6){$\vector(1,0){1}$}

\put(0,3.5){$u_{t-3}$}
\put(2,3.5){$u^1_{t-3}$}
\put(3.9,3.5){$u^{11}_{t-3}$}
\put(0.8,3.6){$\vector(1,0){1}$}
\put(2.8,3.6){$\vector(1,0){1}$}
\put(8,7.5){$u^{2\overbrace{\scriptstyle{1 \cdots 1}}^{t-5}}_1$}
\put(9.1,7.6){$\vector(1,0){1}$}
\put(10.4,7.5){$u^{2\overbrace{\scriptstyle{1 \cdots 1}}^{t-4}}_1$}
\put(2,2.5){$u^2_{t-3}$}
\put(0,1.5){$u_{t-2}$}
\put(0.8,1.6){\vector(1,0){1}}
\put(2,1.5){$u^1_{t-2}$}
\put(0,0.5){$u_{t-1}$}
\put(8,6.5){$u^{2\overbrace{\scriptstyle{1 \cdots 1}}^{t-6}2}_1$}
\put(6.8,7.2){\vector(2,-1){1}}
\put(4.7,7.2){$\vector(2,-1){1}$}
\put(5.9,6.6){$\cdots$}
\put(0.7,3.3){$\vector(2,-1){1}$}

\put(2.2,6.5){\vdots} \put(4,8.8){\vdots}
\put(3.9,8.2){$u^{1(t-3)}_1$} \put(2.5,9.2){$\vector(4,-3){1.2}$}
\put(1.9,5.5){$u^{t-3}_1$} \put(1.8,4.5){$u^{t-2}_1$}
\put(4,5.5){$u^{(t-3)1}_1$} \put(2.8,5.6){$\vector(1,0){1}$}
\put(0.4,9){$\vector(1,-4){1.1}$}
\put(0.5,9.1){$\vector(1,-3){1.1}$}
\put(4.3,9.3){$\vector(2,-1){1}$} \put(5.6,8.7){$\cdots$}
\multiput(5.4,9)(0.2,0.2){3}{\circle*{0.052}}
\multiput(5.4,7)(0.2,0.2){3}{\circle*{0.052}}
\multiput(3.4,7)(0.2,0.2){3}{\circle*{0.052}} \put(0.2,4.5){\vdots}
\put(0.2,4){\vdots} \put(6,-.5){Fig.2} }
\end{picture}

Let $F(i), 1 \le i \le t-1$, denote the connected component of $F_t$
containing $u_i$, hence $F_t=\bigcup\limits_ {i=1}^{t-1}F(i)$.

Without confusing with the above paragraphs we call $u_i$ the $root$
vertex of directed graph $F(i)$ too,\ $1 \le i \le t-1$.

{\bf Lemma 4$^{[7]}$}\quad A weak digraph is an outtree if and only
if it has exactly one vertex with indegree 0, and other vertices
with indegree 1.

The following lemma is an immediate consequence of Lemma 4 and
definition of $F(i)$.

{\bf Lemma 5}\quad $\forall \ i, 1 \le i \le t-1$, $F(i)$ is an
outtree with $root$ $u_i$, then there is a directed path from $u_i$
to any $u^{i_1i_2\cdots i_n}_i$.

Let $F^{'}(i)$ be a disjoint isomorphic graph of $F(i)$ with the
vertex set $V(F^{'}(i))$=\{$v^{i_1i_2\cdots i_n}_i\ :\
i+\sum\limits_{l=1}^n{i_l} \le t-1\}$,\ where $v^{i_1i_2\cdots
i_n}_i$ corresponds to $u^{i_1i_2\cdots i_n}_i \in V(F(i))$.

{\bf Definition 7}\quad For $1 \le i \le t-1$, $F(i)\sqcup F^{'}(i)$
is a digraph with vertex set
$$V(F(i)\sqcup F^{'}(i))=V(F(i)) \cup V(F^{'}(i)),$$
and arc set$$ A(F(i)\sqcup F^{'}(i))=A(F(i))\cup A(F^{'}(i))
\cup\{(u_i,v_i)\}.$$

{\bf Lemma 6}\quad For $1 \le i \le t-2$, $F(i) \cong F(i+1)\sqcup
F^{'}(i+1)$.

{\bf proof}\quad  Let $g$ be a mapping from $V(F(i))$ to
$V(F(i+1)\sqcup F^{'}(i+1))$, i.e.,
$$
g:\ \ V(F(i))\longrightarrow V(F(i+1)\sqcup F^{'}(i+1))
$$
such that
$$
g(x)= \left\{\begin{array}{ll} u_{i+1} ,& x=u_i;\\
\\
u^{(i_1-1)i_2\cdots i_n}_{i+1}, & x=u^{i_1i_2\cdots i_n}_i,i_1 \ge 2;\\
\\
v_{i+1} , & x=u^1_i ;\\
\\
v^{i_2\cdots i_n}_{i+1}, & x=u^{1i_2\cdots i_n}_i.
\end{array}\right.
$$
Then \quad \ $(u_i,u^1_i) \in A(F(i))\Longleftrightarrow
(u_{i+1},v_{i+1}) \in A(F(i+1)\sqcup F^{'}(i+1))$, $$ (u^{1i_2\cdots
i_{n-1}}_i,u^{1i_2\cdots i_n}_i) \in A(F(i))\Longleftrightarrow
(v^{i_2\cdots i_ {n-1}}_{i+1},v^{i_2\cdots i_n}_{i+1}) \in
A(F(i+1)\sqcup F^{'}(i+1)),$$ and for $i_1 \ge 2$, $$
(u^{i_1i_2\cdots i_{n-1}}_i,u^{i_1i_2\cdots i_n}_i) \in A(F(i))
\Longleftrightarrow (u^{(i_1-1) i_2\cdots
i_{n-1}}_{i+1},u^{(i_1-1)i_2\cdots i_n}_{i+1}) \in A(F(i+1)\sqcup
F^{'}(i+1)). $$ Hence $g$ is an isomorphism.\qquad \hfill$\Box$

According to Lemma 6,  it is straightforward to verify
$|V(F(i))|=2^{t-1-i},1 \le i \le t-1$.

{\bf Definition 8}\quad Let $D$ be a digraph, $\{u,v,w\} \subseteq
V(D)$, $(u,v,w)$ is a {\it directed triple } of $D$ if there is a
directed path $P(u,v,w)$ from $u$ to $w$, and contains $v$ as an
inner vertex. Let $S \subseteq V(D)$, $S$ is a 3-{\it cut set} of
$D$ if there is no directed triple of $D$ in $V(D)-S$.

{\bf Lemma 7}\quad $t \ge 3$, $S$ is 3-cut set of $F_t$, then $|S|
\ge 2^{t-3}-1$, and there exists a 3-cut set such that ` = ' holds.

{\bf proof}\quad  If $t=3$, then the length of the longest directed
path of $F_3$ is 1. So the directed triple is $\emptyset$. Hence the
proposition comes true.

Let $t \ge 4$, and suppose that $S$ is a 3-cut set of $F_t$, then $S
\cap V(F(i))=S_i$ is a 3-cut set of $F(i)$, and $S=\bigcup\limits_
{i=1}^{t-1}S_i$. Assume that $S_i$ is the smallest 3-cut set of
$F(i)$, it is easy to see that we need to confirm $|S_i| =
2^{t-i-3}$, $1 \le i \le t-3$.

Let $S_i=\{u^{i_1i_2\cdots i_n}_i: i_1+i_2+\cdots +i_n+i \le
t-3\}$,\ we will prove $S_i$ is the smallest 3-cut set of $F(i)$ by
induction on $i$.

For $i=t-3$, $S_{t-3}=\{u_{t-3}\}$ is the smallest 3-cut set of
$F(t-3)$, the proposition holds. For $i,\ 1 \le i \le t-4$,
hypothesize the conclusion holds for $i+1$. It means that
$S_{i+1}=\{u^{i_1i_2 \cdots i_n}_{i+1}: i_1+i_2+\cdots +i_n+(i+1)
\le t-3\}$ is the smallest 3-cut set of $F(i+1)$.

In the case of $F(i)$, by Lemma 6, $F(i) \cong F(i+1)\sqcup
F^{'}(i+1)$. Since $F(i+1) \cong F^{'}(i+1)$, as the image of
$S_{i+1}$, $S^{'}_{i+1}$ is the smallest 3-cut set of $F^{'}(i+1)$.
So $S_{i+1} \cup S^{'}_{i+1}$ is a 3-cut set of $F(i+1)\sqcup
F^{'}(i+1)$. By the mapping $g$ in Lemma 6, the preimage of $S_{i+1}
\cup S^{'}_{i+1}$ is $\{u^{1i_1i_2\cdots i_n}_i,u^{(i_1+1)i_2\cdots
i_n}_i:1+i_1+i_2+\cdots +i_n+i \le t-3\}=S_i$,\ which is a 3-cut set
of $F(i)$.  $|S_i|=2^{t-i-3}$ follows by induction.

If $S_i$ is not the smallest 3-cut set of $F(i)$, then there is a
3-cut set of $F(i)$, denoted as $T$, such that $|T|=|g(T)|<|S_i|
=2^{t-i-3}$. Since $g(T)$ is a 3-cut set of $F(i+1)\sqcup
F^{'}(i+1)$, and $g(T)=(g(T)\cap V(F(i+1))) \bigcup (g(T)\cap
V(F^{'}(i+1)))$,\ assume $|(g(T)\cap V(F(i+1)))| \le 2^{t-(i+1)-3}-1
< |S_{i+1}|$, contradicting to $S_{i+1}$ is the smallest 3-cut set
of $F(i+1)$.

Hence $S=\bigcup\limits_ {i=0}^{k-1}S_i$ is the smallest 3-cut set
of $F_t$ satisfying $|S|=2^{t-3}-1$.\hfill $\Box$

{\bf Corollary 1}\quad Let $t \ge 4$, $U \subseteq V(F^{\circ}_t)$,
if $|U|>3\cdot2^{t-3}$, then there exists a directed triple in $U$.

{\bf proof}\quad  Assume to the contrary that there is no directed
triple in $U$, then $R_t(M^t(G))-\{u_t\}-U$ is 3-cut set, and
$|R_t(M^t(G))-\{u_t\}-U| < 2^{t-1}-1-3\cdot2^{t-3}=2^{t-3}-1$,
contradicts to $F^{\circ}_t \cong F_t$ and lemma 7. \hfill $\Box$

\section{Main result}

It is straightforward to verify the following lemma by definition of
$M^t(G)$.

{\bf Lemma 8}\quad  Let $v_1,v_2$ be two $t$th $roots$ of
$M^t(K_n)$, and there is a directed path $P(v_1,v_2)$ from $v_1$ to
$v_2$ in $F^{\circ}_t$, then $N(\{v_2,h^{-1}(v_2)\}) \cap T(V(K_n))
\subseteq N(\{v_1,h^{-1}(v_1)\})$, $N(v_2) \cap T(V(K_n)) \subseteq
N(v_1)$.

{\bf Lemma 9}\quad  Let $t \ge 3$, $\chi_c(M^t(K_n))=\frac{k}{2}$,
and a (k,2)-coloring satisfy Lemma 3. If there exists a set of $3$
$t$th $roots$ $\{v_1,v_2,v_3\}$ such that $C(v_i) \subseteq
\delta(v_i) \cap R(M^t(K_n)) \subseteq \{v_i,h^{-1}(v_i)\}$, $1 \le
i \le 3$, then $(v_1,v_2,v_3)$ is not the directed triple of
$F^{\circ}_t$.

{\bf proof}\quad  Suppose on the contrary that $(v_1,v_2,v_3)$ is a
directed triple of $F^{\circ}_t$, and $C(v_i) \subseteq
\{v_i,h^{-1}(v_i)\}$, $v_i \in X_{s_i},1 \le i \le 3$. It is
sufficient to consider two cases as follows:

Case 1.\ $h^{-1}(v_3) \notin C(v_3)$. According to Lemma 3, Lemma 8
and the assumption of this lemma, we get $h^{-1}(v_3) \in
\delta(u_t)$, $N(v_3) \cap \delta(v_2) =\emptyset$. Hence, we can
dye $v_3$ color $s_2$, and obtain a new $(k,2)$-coloring $(Y_0, Y_1,
\cdots, Y_{k-1})$ such that $Y_{s_2}=X_{s_2} \bigcup
\{v_3\}$,$Y_{s_3}=\emptyset $ and $Y_i=X_i,i \neq s_2,s_3$, a
contradiction to Lemma 1.

Case 2.\ $h^{-1}(v_3) \in C(v_3)$. Suppose that $h^{-1}(v_2) \in
C(v_2)$, by Lemma 8 and the assumption of this lemma, we get
$N({v_3,h^{-1}(v_3)}) \cap \delta({v_2,h^{-1}(v_2)}) =\emptyset$.
Hence, we can dye ${v_3, h^{-1}(v_3)}$ color $s_2$ and obtain a new
$(k,2)$-coloring $(Y_0,Y_1,\cdots,Y_{k-1})$ such that
$Y_{s_3}=\emptyset$, a contradiction. So $h^{-1}(v_2) \notin
C(v_2)$, then we can dye $v_2$ color of $v_1$ and also obtain a new
$(k,2)$-coloring $(Y_0,Y_1,\cdots,Y_{k-1})$ such that
$Y_{s_2}=\emptyset$, a contradiction.\hfill $\Box$

{\bf Lemma 10}\quad  Let $G$ be a simple graph with
$V(G)=\{x_1,x_2,\cdots,x_n\}$. The vertex set of $M(G)$ is
$V(M(G))=\{x_1, x_2, \cdots, x_n\}\cup \{x^{'}_1, x^{'}_2, \cdots,
x^{'}_n\}\cup \{u\}$. If $M(G)$ has $(k,d)$-partition, then there is
a $(k,d)$-partition of $M(G)$ such that

1.\ $X_0=\{u\}$;\ \ \

2.\ For some $i$,\ $d \le i \le k-d$,$\forall \ x \in V(G)$, we have
$\{x,x^{'}\}\subseteq X_i$ if $x \in X_i$.

{\bf proof}\quad  For any vertex $v \in V(M(G))-\{u\}$,$v^{'}$ is
the twin of $v$.

By Lemma 3, there exists a $(k,d)$-partition satisfying 2 and $u \in
X_0$. If $x \neq u, x \in X_0,$ and $x^{'} \in X_i, d \le i \le
k-d$, then $x \in V(G)$. We can dye $x$ color $i$. Otherwise, there
exists $y \in \delta(x^{'})$ such that $x^{'}y \notin E(M(G)), xy
\in E(M(G))$. So $y \notin V(G)\cup \{u\}, x^{'}y^{'} \in E(M(G))$.
Hence $y^{'}$ can't be dyed any color, this is impossible. Hence, $x
\notin X_0$, $X_0=\{u_t\}$ and conclusion 1 holds as well.\hfill
$\Box$

{\bf Theorem 1}\quad  Let $t \ge 4,n \ge
\frac{11}{12}2^{t-1}+2t+\frac{1}{3}$,\ then $\chi_c(M^t(K_n))=
\chi(M^t(K_n))=n+t$.

{\bf proof}\quad  Suppose that $\chi_c(M^t(K_n))=\frac{k}{d}, d \neq
1$, by Lemma 2
\begin{eqnarray*}
d-1 &\le & \frac{2^t-2}{n-3}\le
\frac{2^t-2}{\frac{11}{12}2^{t-1}+2t-\frac{8}{3}}
\le 3+\frac{2^t-2-3(\frac{11}{12}2^{t-1}+2t-\frac{8}{3})}{\frac{11}{12}2^{t-1}+2t-\frac{8}{3}}\\
&=&3-\frac{3\cdot 2^{t-3}+6t-6}{\frac{11}{12}2^{t-1}+2t-\frac83}<3
\end{eqnarray*}
holds for $t \ge 3$. So $d \le 3$.

If $d=3$, there exists a $(3n+3t-i,3)$-coloring of $M^t(K_n)$, $i=1$
or $2$. According to box principle, there are at least
$5n-(3n+3t-i)=2n-3t+i$ color classes in the intersection sets of
$d$-fields of different initial vertices. By Lemma 10, there are at
most 3 such color classes in $\delta(u_t)$. So
$V(M^t(K_n))-\delta(u_t)$ have at least $2n-3t+i-3$ such color
classes in which contains at least one vertex in
$R_t(M^t(G))-\{u_t\}$ respectively. Thus $2n-3t+i-3 \le 2^{t-1}-1$,
i.e.,

$$n \le 2^{t-2}+\frac{3}{2}t+\frac{2-i}{2}\le2^{t-2}+\frac{3}{2}t+\frac{1}{2}.$$
If $t \ge
4$,then$$\frac{11}{12}2^{t-1}+2t+\frac{1}{3}-(2^{t-2}+\frac{3}{2}t+\frac{1}{2})=\frac{5}
{6}2^{t-2}+\frac{1}{2}t-\frac{1}{6}>0.$$

So$$ n \ge \frac{11}{12}2^{t-1}+2t+\frac{1}{3}
>2^{t-2}+\frac{3}{2}t+\frac{1}{2} \ge n,
$$a contradiction.

Suppose that $d=2$, then $M^t(K_n)$ has a $(2n+2t-1,2)$-coloring
$(X_0,X_1,\cdots,X_{2n+2t-2})$ satisfying Lemma 10. According to box
principle, there are at least $3n-(2n+2t-1)=n-2t+1$ color classes in
the intersection set of $d$-fields of different initial vertices. By
Lemma 10, there are at most 1 such color class in $\delta(u_t)$. By
definition 5, any of these color classes we get is in $R(M^t(K_n))$,
and has at least one $t$th root.

Now we discuss such color classes  in $V(M^t(K_n))-\delta(u_t)$. It
is easy to see that the number of these color classes is at least
$n-2t$.\ Suppose that the number of these color classes that contain
exactly 1 $t$th root is $y$, the number of those which contain at
least 2 $t$th roots is $x$ and  the number of $t$th roots which is
different from $\{ u_t\} $ and not in such intersection sets is $z$,
then we have
$$x+y \ge n-2t,\ \ \ \ 2x+y+z \le 2^{t-1}-1.$$
Hence
$$y \ge 2n-2^{t-1}-4t+1+z,\ \ \ z \le 2^{t-1}+2t-n-1.$$

Since $d=2,$ the color classes in the intersection set of $d$-field
of different initial vertices are mutually disjoint. Hence, in the
set of $y$ classes which contain exactly 1 $t$th root,\ the number
of classes whose $d$-fields contain other $t$th roots is at most
$2z$,\ furthermore there are at most 2 classes whose $d$-fields
contain $r$th roots, $1 \le r \le t-1$, $X_2$ and $X_{2n+2t-3}$.
Hence, there are at least $y-2z-2$ color classes whose $d$-fields
contain exactly 1 $t$th root.

By Corollary 1, if $y-2z-2 > 3\cdot2^{t-3}$, then there exists a
directed triple of $F^{\circ}_t$ in the set of $t$th roots that
contained in these $y-2z-2$ color classes, contradict to Lemma 9.
Hence, $y-2z-2 \le 3\cdot2^{t-3}$, we have
\begin{eqnarray*}
&     &(2n-2^{t-1}-4t+1+z)-2z-2 \le y-2z-2 \le 3\cdot2^{t-3}\\
&\Longrightarrow & 2n-2^{t-1}-4t-1-z \le 3\cdot2^{t-3}\\
&\Longrightarrow & 2n-2^{t-1}-4t-1-(2^{t-1}+2t-1-n) \le 3\cdot2^{t-3}\\
&\Longrightarrow & 3n \le 11\cdot2^{t-3}+6t\\
&\Longrightarrow & n \le \frac{11}{12}\cdot2^{t-1}+2t,
\end{eqnarray*}
contradict to $n \ge \frac{11}{12}2^{t-1}+2t+\frac{1}{3}$. So $d
\neq 2$.

In other words, if $n \ge \frac{11}{12}2^{t-1}+2t+\frac{1}{3}$, then
$d=1$, $\chi_c(M^t(K_n))=\chi(M^t(K_n))=n+t$. \qquad $\Box$

\end{document}